\documentclass[12pt]{article}
\usepackage{theorem,amsfonts, amssymb}
\newtheorem{theorem}{Theorem}[section]
\newtheorem{proposition}[theorem]{Proposition}
\newtheorem{lemma}[theorem]{Lemma}
\newtheorem{corollary}[theorem]{Corollary}
\theorembodyfont{\upshape}
\newtheorem{definition}{Definition}[theorem]
\newtheorem{example}[theorem]{Example}
\newtheorem{remark}[theorem]{Remark}
\newtheorem{proof}{Proof}

\newtheorem{acknowledgement}{Acknowledgement}

\newcommand{\bt}{\begin{theorem}}
\newcommand{\et}{\end{theorem}}
\newcommand{\bl}{\begin{lemma}}

\newcommand{\el}{\end{lemma}}
\newcommand{\bp}{\begin{proposition}}
\newcommand{\ep}{\end{proposition}}
\newcommand{\bd}{\begin{definition}}
\newcommand{\bex}{\begin{example}}
\newcommand{\eex}{\end{example}}
\newcommand{\ed}{\end{definition}}
\newcommand{\br}{\begin{remark}}
\newcommand{\er}{\end{remark}}
\newcommand{\bc}{\begin{corollary}}
\newcommand{\ec}{\end{corollary}}
\newcommand{\bo}{\begin{proof}}
\newcommand{\eo}{\end{proof}}
\newcommand{\be}{\begin{enumerate}}
\newcommand{\ee}{\end{enumerate}}

\newcommand{\Aut}{{\rm Aut}}
\newcommand\du{\,{\rm d}}

\newcommand{\supp}{{\rm supp}}
\newcommand{\Z}{{\mathbb Z}}

\newcommand{\N}{{\mathbb N}}

\newcommand{\K}{{\cal K}}
\newcommand{\A}{{\cal A}}

\title{Orbits of Distal actions on Locally Compact Groups}
\author{Riddhi Shah}
\date{October 22, 2010} 

\begin{document}
\maketitle

\let\epsi=\epsilon
\let\vepsi=\varepsilon
\let\lam=\lambda
\let\Lam=\Lambda
\let\ap=\alpha
\let\vp=\varphi
\let\ra=\rightarrow
\let\Ra=\Rightarrow
\let\da=\downarrow
\let\Llra=\Longleftrightarrow
\let\Lla=\Longleftarrow
\let\lra=\longrightarrow
\let\Lra=\Longrightarrow
\let\ba=\beta
\let\ga=\gamma
\let\Ga=\Gamma
\let\un=\upsilon
\let\mi=\setminus
\let\ol=\overline
\let\ot=\odot

\begin{abstract}
We discuss properties of orbits of (semi)group actions on locally compact 
groups $G$, In particular, we show that if a compactly generated locally compact abelian 
group acts distally on $G$ then the closure of each of its orbits is a minimal closed invariant set
(i.e.\ the action has [MOC]). We also show that for such an action distality is preserved
if we go modulo any closed normal invariant subgroup and hence [MOC] is also preserved. 
We also show that any semigroup action on $G$ has [MOC] if and only if  the 
corresponding actions on a compact invariant metrizable subgroup $K$ and on the quotient 
space $G/K$ has [MOC].

\end{abstract}

\begin{section}{Introduction}

Let $X$ be a Hausdorff space and $\Ga$ be a (topological) semigroup
acting continuously on $X$ by continuous self-maps. 
The action of $\Ga$ on $X$ is said to be {\it distal} if for any two distinct
points $x, y \in X$, the closure  of $\{ (\ga (x), \ga (y)) \mid \ga
\in \Ga \}$ does not intersect the diagonal $\{ (a, a ) \mid a \in X \}$;
it is said to be {\it pointwise distal} if
for each $\gamma\in\Gamma$, the action of $\{\gamma ^n\}_{n\in\N}$ on $X$ is
distal.  The $\Gamma$-action on $X$ is said to have [MOC] 
(minimal orbit closures) if the closure of every $\Gamma$-orbit is a minimal 
closed $\Gamma$-invariant set, i.e.\ for $x,y\in X$, if $y\in\ol{\Ga(x)}$ then 
$\ol{\Ga(y)}=\ol{\Gamma(x)}$. The notion of distality was introduced by Hilbert
(cf.\ Ellis \cite{El}, Moore \cite{Mo}) and studied by many in different contexts, 
(see Abels \cite{Ab1}-\cite{Ab2},  Furstenberg \cite{Fu}, 
Raja-Shah \cite{RS} and the references cited therein).

Let $G$ be a locally compact (Hausdorff) group and let $e$ denote the
identity of $G$.  Let $\Ga$ be a semigroup acting continuously on $G$
by endomorphisms.  Then $\Gamma$-action on $G$ is distal if and only if
$e\not\in\overline {\Gamma x}$ for all $x\in G\mi\{e\}$.  Note that if $\Ga$-action
on $G$ has [MOC], then it is distal; for if $e\in\ol{\Gamma x}$, then $\{e\}=\ol{\Gamma e}
=\ol{\Gamma x}$ and hence $x=e$. What we are interested in is the converse. 
If $\Gamma$-action on $G$ is distal, does it have [MOC]? The answer is known to be affirmative in any
of the following cases: (1)  $G$ is compact (2)  $\Gamma$ is compact, (3) $G$ is a connected Lie group and 
$\Gamma$ is a subgroup of $\Aut(G)$ (4) G is discrete, or more generally, all $\Ga$-orbits are closed. 
If $\Ga$ is a group and if $\Ga'$ is a closed co-compact normal subgroup, then 
$\Ga$-action on $G$ has [MOC] if and only if $\Ga'$-action on $G$ has [MOC] (cf.~\cite{Mo}); it is 
easy to see that the same equivalence is true for distality. For a locally compact group $G$ and a group 
$\Gamma\subset\Aut(G)$ acting distally on $G$, the  answer to the above question is not known. But
in case of a certain kind of $\Gamma$, we get the following: 

\bt\label{main}
Let $G$ be a locally compact group and let $\Gamma$ be a compactly generated locally compact 
abelian group such that $\Gamma$ acts on $G$ by automorphisms. Then the following are equivalent: 
\begin{enumerate}
\item The $\Gamma$-action on $G$ is distal
\item The $\Gamma$-action on $G$ has {\rm [MOC]}.
 \end{enumerate}
\et

Let us now discuss general actions on compact spaces.
For a compact space $K$, let $\Gamma$ be a semigroup of continuous bijective self-maps of
$K$. Then $\Gamma$ is a subsemigroup of $C(K)$, the group of all continuous 
bijective self-maps on $K$. Let $[\Gamma]$ be the group generated by $\Gamma$ in $C(K)$. 
We know that $\Gamma$-acts distally on $K$ if and only if $E(\Gamma)$, the 
closure of $\Gamma$ in $K^K$ with weak topology, is a group (cf.\ \cite{El}); it is obviously compact 
since $K^K$ is so. Then $E(\Gamma)=E([\Gamma])$. Moreover, 
for any $x\in K$, $\ol{\Gamma(x)}=E(\Gamma)(x)=E([\Gamma])(x)$. So for a compact space $K$ 
and $\Gamma$ and $[\Gamma]$ as above, the following are equivalent: 
\begin{enumerate} 
\item $\Gamma$-action on $K$ is distal.
\item $[\Gamma]$-action on $K$ is distal.
\item $\Gamma$-action on $K$ has [MOC].
\item $[\Gamma]$-action on $K$ has [MOC].
\end{enumerate}

In particular, if $G$ is a locally compact group and $\Gamma$ a semigroup in $\Aut(G)$ 
such that $\Gamma$ keeps a closed co-compact subgroup $H$ of $G$ invariant 
(i.e.\ $\gamma(H)=H$ for all $\gamma\in\Gamma$), then the above equivalence is also 
true for the actions of $\Gamma$ and $[\Gamma]$ on $G/H$. Note that for any $\Gamma$-action on
$G$, the coresponding $\Gamma$-action on the homogeneous space $G/H=\{xH\mid x\in G\}$ is 
canonically defined as 
$\gamma(xH)=\gamma(x)H$ for all $\gamma\in\Ga$; it is well-defined since $H$ is $\Gamma$-invariant. 

In \cite{RS}, it is shown that distality of a semigroup action is preserved by factor actions modulo  
compact invariant subgroups. We show that a similar result holds for [MOC], (see also Remark \ref{rem1}).

\bt\label{cpt}
Let $G$ be a locally compact group and let $\Gamma$ be a 
subsemigroup of $\Aut(G)$. Let $K$ be a compact metrizable $\Gamma$-invariant 
subgroup of $G$. Then $\Gamma$-action on $G$ has {\rm [MOC]} if and only if
$\Gamma$-action on both $K$ and $G/K$ has {\rm [MOC]}.
\et

The following result is about factor actions modulo closed normal invariant subgroups.

\bt\label{normal}
Let $G$ and $\Gamma$ be as in Theorem~\ref{main}. Let $H$ be a closed normal $\Gamma$-invariant 
subgroup of $G$. Then $\Gamma$-action on $G$ has {\rm [MOC]} if and only if $\Gamma$-action on both 
$H$ and $G/H$ has {\rm [MOC]}.
\et

We will later show that a similar result holds for distality for a larger class of $\Gamma$.

A locally compact  group $G$ is said to be {\it distal} (resp.\
{\it pointwise distal}) if the conjugacy action of $G$ on $G$ is distal
(resp.\ pointwise distal). A distal group is obviously pointwise distal. 
It can easily be seen that the class of distal groups is closed under
compact extensions.  Abelian groups, discrete groups and
compact groups are obviously distal.
Nilpotent groups, connected groups of
polynomial growth are distal (cf.\ \cite{Ro}) and p-adic Lie
groups of type $R$ and p-adic Lie groups of polynomial growth are
pointwise distal (cf.\ Raja \cite{R1} and \cite{R2}).

In \cite{RS}, we have shown that any locally compact
group is pointwise distal if and only if it has shifted convolution property; i.e.\ 
for any probability measure $\mu$ on $G$, whose concentration functions
do not converge to zero, there exists $x\in\supp\mu$, the support of $\mu$, such that 
$\mu^nx^{-n}\to\omega_H$, the
Haar measure of some compact group $H$ which is normalised by $\supp\mu$. 
For a probability measure $\mu$ on $G$, the $n$-th convolution
function of $\mu$ is defined as $f_n(\mu, C)=\sup_{g\in G}\mu^n(Cg)$, for any compact subset $C$
of $G$. We say that the concentration functions of $\mu$ do not converge to zero if there exists a compact set
$C$ such that $f_n(\mu, C)\not\to 0$ as $n\to\infty$, (see \cite{RS} for more details).
The following corollary is a consequence of Theorem 6.1 of \cite{RS} and Theorem \ref{main}. 

\bc\label{distal} Let $G$ be a locally compact group. 
Then the following are equivalent:
\begin{enumerate} 
\item $G$ is pointwise distal. 
\item $G$ has shifted convolution property.
\item For every $g\in G$, the conjugation action of $\{g^n\}_{n\in\Z}$ on $G$ has {\rm [MOC]}. 
\end{enumerate}
\ec

A locally compact group $G$ is said to be a {\it generalised} $FC^-$-{\it group} (resp.\  $FC^-$-{\it nilpotent}) if
$G$ has closed normal subgroups $\{G=G_0,\ldots, G_n=\{e\}\}$ such that  $G_{i+1}\subset G_i$ and 
$G_i/G_{i+1}$ is a compactly generated group with relatively compact conjugacy classes (resp.\ 
every orbit of the conjugacy action of $G$ on $G_i/G_{i+1}$ is relatively compact) for all 
$i=0,1,\ldots,n-1$. Any compactly generated group $G$ has polynomial growth if and only if it is 
$FC^-$-nilpotent; and it is a generalised $FC^-$-group (cf.\ \cite{L2}). Any compactly generated 
abelian group (resp.\ any polycyclic group) is a generalised $FC^-$-group. More generally, any
compactly generated group with polynomial growth is a generalised $FC^-$-group. Note that 
generalised $FC^-$-groups are compactly generated (cf.\ \cite{L2}, Proposition 2).

Recall that a subgroup $\Ga$ of $\Aut(G)$ is said to be equicontinuous (at $e$) if and only if there exists a 
neighbourhood base at $e$ consisting of $\Ga$-invariant neighbourhoods; in case of totally disconnected
groups, this is equivalent to the existence of a neighbourhood base at 
$e$ consisting of compact open $\Ga$-invariant subgroups. If $\Ga$ is compact, then it is easy to 
see that $\Ga$ is equicontinuous. If $G$ is a totally disconnected group and if $\Ga$ has
a polycyclic subgroup of finite index and it acts distally on $G$, then $\Ga$ is
equicontinuous (cf.\ \cite{JR}, Corollary 2.4). If any group $\Ga$ acts on $G$ by automorphisms and 
its image in $\Aut(G)$ is equicontinuous then we say that $\Ga$-action on $G$ is equicontinuous.

For a totally disconnected locally compact group $G$, 
we have the following result:

\bp\label{tds} Let $G$ be a totally disconnected locally compact group and let 
$\Gamma$ be a generalised $FC^-$-group which acts on $G$ by automorphisms. Then 
the following are equivalent.
\begin{enumerate}
\item $\Gamma$-action on $G$ is distal.
\item $\Gamma$-action on $G$ has {\rm [MOC]}.
\item $\Gamma$-action on $G$ is equicontinuous.
 \end{enumerate}
\ep

In Section 2, we discuss factor actions modulo compact (resp.\ closed normal) invariant groups and prove Theorem \ref{cpt}, 
Proposition \ref{tds} and an analogue of Theorem \ref{normal} for distal actions of a more general class of groups. 
In Section 3, we prove the equivalence of distality and [MOC] of certain actions, namely, Theorem \ref{main}. 
Note that if $\Gamma$ acts on $G$ by automorphisms, for convenience, $\Ga$ is often equated with 
its image in $\Aut(G)$, whenever there is no loss of any generality. 

\end{section}

\begin{section}{Orbits of Factor Actions}

In this section we discuss [MOC] of factor actions modulo compact invariant groups and modulo closed
normal invariant groups. We first show that [MOC] is preserved if we go modulo a compact invariant subgroup
by proving Theorem \ref{cpt}. Before that we prove a proposition which proves
a special case of the theorem in case the compact subgroup is a Lie group. 

\bp\label{cptl}
Let $G$ be a locally compact group and let $\Gamma$ be a 
subsemigroup of $\Aut(G)$. Let $K$ and $L$ be compact $\Gamma$-invariant subgroups
of $G$ such that $L$ is a normal subgroup of $K$ and $K/L$ is a Lie group. Then 
$\Gamma$-action on $G/L$ has {\rm [MOC]} if and only if
$\Gamma$-action on both $G/K$ and $K/L$ has {\rm [MOC]}.
\ep

\bo {\bf Step 1} Let $G$, $\Ga$, $K$ and $L$ be as in the hypothesis. 
One way implication ``only if'' is easy to prove.
Suppose $\Gamma$-action on $G/L$ has [MOC]. Then clearly $\Gamma$-action on
$K/L$ also has [MOC], as $K$ is closed and $\Gamma$-invariant. 
Now we want to show that $\Gamma$-action on $G/K$ has [MOC]. 
Let  $x\in G$ and let $yK\in\ol{\Gamma(xK)}$ in $G/K$ for some $y\in G$. 
Then $yK\subset\ol{\Gamma(x)K}=\ol{\Gamma(x)}K$ and hence
$yk\in\ol{\Gamma(x)}$ for some $k\in K$. In particular, we get that $ykL\subset\ol{\Gamma(x)}L=
\ol{\Gamma(x)L}$ as $L$ is compact. Hence $ykL\in\ol{\Gamma(xL)}$
in $G/L$. Since $\Gamma$-action on $G/L$ has [MOC], we get that $\ol{\Gamma(xL)}=\ol{\Gamma(ykL)}$
and hence $x\in \ol{\Gamma(y)K}$ as $k\in K$, $L\subset K$
and both $L$ and $K$ are $\Gamma$-invariant. This implies that $xK\in\ol{\Gamma(yK)}$ in
$G/K$ and hence $\Gamma$-action on $G/K$ has [MOC]. Note that the condition that $K/L$ is a Lie group
is not used in the proof of the ``only if'' statement.

\medskip
\noindent{\bf Step 2} Now we prove the ``if'' statement. Suppose $\Ga$-action on both $G/K$ and $K/L$ has [MOC].
This implies that $\Ga$-action on both $G/K$ and $K/L$ is distal and hence $\Ga$-action
on $G/L$ is distal; (this is easy to see from the proof of Theorem 3.1 in \cite{RS}). 

For any $g\in G$, let $g'=gL$. The map $g\mapsto g'$ is a continuous proper map from $G$
to $G/L$. Let $x\in G$ and let $y'\in\ol{\Gamma(x')}$ for some $y\in G$. We want to show that 
$x'\in\ol{\Ga(y')}$. Then $yK\in\ol{\Ga(xK)}$, and as $\Ga$-action on $G/K$ has [MOC], $xK\in\ol{\Ga(yK)}$.
This implies that $xk\in \ol{\Ga(y)}$ for some $k\in K$, and hence, $x'k'\in\ol{\Ga(y')}$. Let $\{\gamma_d\}$ and $\{\beta_d\}$ be nets
in $\Ga$ such that $\gamma_d(x')\to y'$ and $\beta_d(y')\to x'k'$. 

\medskip
\noindent{\bf Step 3} Let $\Ga_0$ be the closure of image of $\Ga$ in $\Aut(K/L)$. Suppose $\Ga_0$ is compact. 
Then $\Ga_0$, being a compact semigroup, is a group. Let $\beta$ and $\gamma$ be limit points
of images of $\{\beta_d\}$ and $\{\gamma_d\}$ in $\Ga_0$ respectively. Then 
$$
\gamma_d(x'k')\to y'\gamma(k')\in \ol{\Ga(y')}\ \ \mbox{and}\ \ \beta_d(y'\ga(k'))\to x'k'\alpha(k')\in\ol{\Ga(y')},$$
where $\alpha=\beta\gamma\in\Aut(G/K)$.  
Similarly we get that for 
$$
k_n=k'\alpha(k')\cdots\alpha^{n-1}(k')\in K/L,\ \  x_n=x'k_n\in\ol{\Ga(y')},\ \ \mbox {
for all } n\in\N.$$ 
As $\Ga_0$ is a compact group, there exists a sequence  $\{n_j\}\subset\N$ 
such that $\alpha^{n_j}\to I$, the identity of $\Aut(K/L)$. 
Passing to a subsequence if necessary, we may assume that
$k_{n_j}\to c'=cL\in K/L$, for some $c\in K$. Hence $x'c'\in\ol{\Ga(y')}$. 
Now as $\alpha^{n_j}\to I$, 
$$k_{2n_j}=k_{n_j}\alpha^{n_j}(k_{n_j})\to (cL)^2=c^2L.$$
Similarly, for all $m\in\N$, 
$$k_{mn_j}=k_{n_j}\alpha^{n_j}(k_{n_j})\cdots\alpha^{(m-1)n_j}(k_{n_j})\to c^mL\in K/L$$ 
and $xc^mL\in\ol{\Ga(yL)}$. Since $K/L$ is a compact (Lie) group, 
$e'=eL$ is in the closure of $\{c^mL\}_{m\in\N}$ in $K/L$ and hence 
$x'\in\ol{\Ga(y')}$, i.e.\ $\ol{\Ga(x')}=\ol{\Ga(y')}$. Hence $\Ga$-action on $G/L$ has [MOC].

In particular, since $K^0L/L$ is the connected component of $K/L$, $K/K^0L$ is finite, and hence, $\Aut(K/K^0L)$
is finite. Arguing as above for $K^0L$ in place of $L$, we get that $G/K^0L$ has [MOC] and we may assume that 
$K=K^0L$, i.e.\ $K/L$ is connected.

\medskip
\noindent{\bf Step 4} Now Let $Z$ be the subgroup of $K$ such that $L\subset Z$ and $Z/L$ is the center of $K/L$. Then $Z$
and $Z^0L$ are closed and $\Ga$-invariant. Moreover, $K/Z$ is a connected semisimple Lie group and hence its 
automorohism group is compact. Therefore arguing as in Step 3 for $Z$ in place of $L$, we get that $\Ga$-action on 
$G/Z$ has [MOC], and since $Z/Z^0L$ is finite, $\Ga$-action on $G/Z^0L$ also has [MOC]. 
Now replacing $K$ by $Z^0L$, we may assume that $K/L$ is a connected abelian Lie group.

Let $[\Ga]$ be  the group generated by $\Ga$ in $\Aut(K/L)$. Then $[\Ga]$ also acts distally
on $K/L$. By Lemma 2.5 of \cite{Ab2}, 
there exists a finite set of compact (normal) $[\Ga]$-invariant  subgroups $\{K_0,\ldots,K_n\}$ in $K$ such 
that $K=K_0\supset K_1\supset\cdots\supset K_n=L$ and the image of $[\Ga]$ in $\Aut(K_i/K_{i+1})$ is finite
for each $i\in\{0,\ldots, n-1\}$. 
Arguing as in Step 3 for $K_1$ in place of $L$, we get that $\Ga$-action on $G/K_1$, has [MOC]. 
Since the image of $\Ga$ in $\Aut(K_i/K_{i+1})$ is finite, using the above argument repeatedly for $K_i/K_{i+1}$ 
in place of $K/L$, we get that $\Ga$-action on $G/K_{i+1}$ has [MOC], $1\leq i\leq n-1$. Since $K_n=L$, we have 
that $\Ga$-action on $G/L$ has [MOC]. \hfill{$\Box$}
\eo

\bo $\!\!\!\!\!$ {\bf of Theorem \ref{cpt}} Let $G$, $\Ga$ and $K$ be as in the hypothesis.
As in the proof of Proposition \ref{cptl}, the ``only if'' statement is obvious.   
Now we prove the ``if'' statement. Suppose that $\Ga$-action on $G/K$ and $K$ has [MOC].
Hence $\Ga$-actions on $G/K$, $K$ and $G$ are distal. Let ${\cal K}$ consists of 
closed (compact) $\Gamma$-invariant subgroups 
$C$ of $K$ such that $\Gamma$-action on $G/C$ has [MOC]. Then ${\cal K}$ is nonempty as 
$K$ belongs to ${\cal K}$. We put  an order on $\K$ by set inclusion. 
Let $\A=\{K_d\}$ be a totally ordered subset of $\K$. We show that $K'=\cap K_d\in \K$. 

For any $x\in G$ and $y\in\ol{\Ga(x)}K'$, we show that
$\ol{\Gamma(x)}K'=\ol{\Gamma(y)}K'$. We know that $\ol{\Ga(x)}K_d=\ol{\Ga(y)}K_d$ for each $d$. 
First we show that $\cap_d\ol{\Gamma(x)}K_d=\ol{\Gamma(x)}K'$. One way inclusion is obvious. 
Let $a\in \cap_d\ol{\Gamma(x)} K_d$. Then $C_d=\ol{\Ga(x)}\cap aK_d\ne
\emptyset$ for all $d$.  Here, $\A'=\{C_d\}$ is a collection of compact sets and intersection of finitely many
subsets in $\A'$ is nonempty since $\A$ is totally ordered. Hence $\cap_d C_d$ is nonempty. But
$$
\cap_d C_d=\cap_d (\ol{\Ga(x)}\cap aK_d)=\ol{\Ga(x)}\cap(\cap_d aK_d)=\ol{\Ga(x)}\cap aK'\ne\emptyset.$$
Hence $a\in\ol{\Ga(x)}K'$. Therefore, $\cap_d\ol{\Gamma(x)}K_d=\ol{\Gamma(x)}K'$. 
Similarly, $\cap_d\ol{\Gamma(y)}K_d=\ol{\Gamma(y)}K'$. This implies that $\ol{\Gamma(x)}K'=\ol{\Gamma(y)}K'$ 
and hence $\Gamma$-action on $G/K'$ has [MOC], i.e.\ $K'\in\K$.

By Zorn's Lemma, there exists a minimal element in $\K$, say $M$. 
Here, $M$ is a compact $\Gamma$-invariant subgroup of $K$ such that $\Ga$-action on
$G/M$ has [MOC] and there is no proper subgroup of $M$ in $\K$.  We show that 
$M=\{e\}$. If possible suppose $M$ is nontrivial.  Since $M\subset K$ is compact and metrizable and since 
$\Gamma$-action on $M$ is distal,  it is not ergodic and there exists a (nontrivial) irreducible 
unitary representation $\chi$ of $M$ such that $\chi\Gamma$ is finite upto 
equivalence classes (cf.\  \cite{Be}, Theorem 2.1, see also \cite{R3} as the action of the group $[\Ga]$ generated by $\Ga$ 
is also distal). 
Let  $L=\cap_{\gamma\in\Ga} \ker(\chi\gamma)$ Then $L$ is a proper closed (compact) normal $\Gamma$-invariant 
subgroup of $M$ and since $\chi\Gamma$ is finite upto equivalence classes, $M/L$ is a (compact) 
Lie group. Moreover, $\Gamma$-action on $M/L$ is distal (cf.\ \cite{RS}, Theorem 3.1) and hence it has 
[MOC]. By Proposition \ref{cptl}, we get that  
$\Gamma$-action on $G/L$ has  [MOC]. Hence $L\in {\cal K}$, a contradiction
to the minimality of $M$ in ${\cal K}$. Hence $M=\{e\}$ and $\Ga$-action on $G$ has [MOC]. \hfill{$\Box$}
\eo

\br \label{rem1}
1. In Theorem \ref{cpt}, if $G$ is first countable then $K$ is also first countable and hence metrizable.

\smallskip
\noindent 2. Theorem \ref{cpt} holds in case $\Gamma$ is a locally compact $\sigma$-compact group, 
(for e.g.\ $\Gamma=\Z$) and $K$ is not (necessarily) metrizable. As in this case, the group $M$ as above is not
necessarily metrizable. Here, 
$\Gamma\ltimes M$ is locally compact and $\sigma$-compact and hence $M$ has arbitrarily small compact normal 
$\Gamma$-invariant subgroups $M_d$ such that $\cap_d M_d=\{e\}$ and $M/M_d$ is second countable and hence metrizable  
(cf.\ \cite{HR}, Theorem 8.7). Now from 
Theorem 3.1 of \cite{RS}, if $\Gamma$-action on $M$ is distal then the corresponding $\Gamma$-action on $M/M_d$ is
also  distal and hence not ergodic and we get a proper closed normal $\Gamma$-invariant subgroup (of $M/M_d$, and hence,) 
of $M$, denote it by $L$ again, such that $M/L$ is a Lie group. Now the assertion is obvious from the above proof.
Note that any compactly generated locally compact group is $\sigma$-compact.
\er

The following corollary follows from Theorem 3.1 in \cite{RS}, Theorem 1.1 in \cite{Ab2} and Theorem \ref{cpt} above
since every connected locally compact group has a unique maximal compact normal (characteristic) subgroup
such that the quotient is a connected Lie group.

\bc\label{mainc1} Let $G$ be a connected locally compact first countable group. Let $\Gamma$ be a
subgroup of $\Aut(G)$. Then $\Ga$-action on $G$ is distal if and only if it has {\rm [MOC]}.
\ec

We now show that [MOC] is preserved by factors modulo closed normal invariant group. Before that we prove Proposition
\ref{tds} and a Lemma which will be useful in proving Theorem \ref{normal-d} below and also Theorem \ref{main}.

\bo {\bf of Proposition \ref{tds}} Let $G$ be a locally compact totally disconnected group and let $\Ga$ be a generalised 
$FC^-$-group acting on $G$ by 
automorphisms. Let $\Ga_0=\{\ga\in \Ga\mid \ga(x)=x \mbox{ for all }x\in G\}$. Then $\Ga_0$ is a closed normal 
subgroup of $\Ga$, $\Ga/\Ga_0$ is isomorphic to a subgroup of $\Aut(G)$. Also, $\Ga/\Ga_0$ is a 
generalised $FC^-$-group. It is easy to see that we can replace $\Ga$ by $\Ga/\Ga_0$
and assume that $\Ga\subset\Aut(G)$. We prove that $(1)\Rightarrow (3)\Rightarrow (2)\Rightarrow (1)$. 

Suppose $\Ga$ acts distally on $G$. As $\Ga$ is 
totally disconnected, it has a compact open normal subgroup $C$ such that $\Ga/C$ has a 
polycyclic subgroup of finite index (cf.\ \cite{L2}). Since $C$ is equicontinuous, By Lemma 2.3 of \cite{JR},   
$\Gamma$-action on $G$ is also equicontinuous, (see also `Note added in Proof' in \cite{JR} for 
non-metrizable groups). Now $G$ has a neighbourhood base 
at $e$ consisting of open compact subgroups $K_d$ which are $\Ga$-invariant  
and $\cap_d K_d=\{e\}$. Since $G/K_d$ is discrete, $\Ga$-action on $G/K_d$ has [MOC]. 
Let $x\in G$ and $y\in\ol{\Ga(x)}$. Then $\ol{\Ga(x)}K_d=\Ga(x)K_d=\Ga(y)K_d=\ol{\Ga(y)}K_d$ 
as $K_d$ is open for all $d$.  $\ol{\Ga(x)}=\cap_d \Ga(x)K_d=\cap_d \Ga(y)K_d=\ol{\Ga(y)}$. 
This proves that $\Ga$-action on $G$ has [MOC]. We know that [MOC] implies distality. \hfill{$\Box$}
\eo

\bl \label{lem}
Let $G$ be a locally compact group. Let $\Gamma$ be a group acting on $G$ by automorphisms. 
Suppose that $\Ga$-action
on $G/G^0$ is equicontinuous. Then there exist open (resp.\ compact) $\Gamma$-invariant 
subgroups $H_d$ (resp.\ $K_d$) such that $H_d=K_d G^0$, $K_d$ is the maximal
compact normal subgroup of $H_d$, $K_d\cap G^0=\cap_d K_d$ is the maximal compact 
normal $\Ga$-invariant subgroup of $G^0$. In particular, if $G^0$ has no nontrivial compact 
normal subgroup, then $K_d$ is totally disconnected and $H_d=K_d\times G^0$, a direct product, 
for all $d$. 
\el

\bo
Since $\Ga$-action on $G/G^0$ is equicontinuous, there exist  open almost connected
$\Gamma$-invariant subgroups $H_d$ such that $\{H_d/G^0\}$ form a neighbourhood base 
at the identity in $G/G^0$ consisting of compact open subgroups. 

Choose $H=H_d$ for some fixed $d$. Since $H$ is Lie projective, there exists a compact 
normal subgroup $C$ in $H$ such that  $H/C$ is a Lie group with finitely many connected 
components.  Hence $H$ has a maximal compact normal subgroup, we denote it by $C$ again.
Then $C$ is characteristic in $H$, and in particular, it is $\Ga$-invariant. Let $H'=CG^0$. It is
an open $\Ga$-invariant subgroup in $G$ and $K=C\cap G^0$ is the maximal compact 
normal subgroup of $G^0$.  Since $H'/G^0$ is compact and open in $G/G^0$, passing to a 
subnet, we may assume that $H_d\subset H'$ for all $d$.  Let $K_d=C\cap H_d$. 
Then $K_d$ is a compact normal $\Gamma$-invariant subgroup in $H_d$ 
and $H_d=K_d G^0$ as $G^0\subset H_d$. As $K=C\cap G^0\subset H_d$,  
$K=K_d\cap G^0$ and $K_d$ is the maximal compact normal 
subgroup in $H_d$ for every $d$. Also, since $\cap_d H_d=G^0$, we get that $\cap_d K_d=K$.
Moreover, if $K_d\cap G^0=K$ is trivial, then $K_d$ is totally disconnected and 
$H_d=K_d\times G^0$ as both $K_d$ and $G^0$ are normal in $H_d$, for all $d$. \hfill{$\Box$}
\eo

\smallskip
 To prove Theorem \ref{normal}, in view of Theorem \ref{main}, it is enough if we prove the same statement for
 distal actions. Here, we prove the following for  distal actions of a more general class of groups. 

\bt\label{normal-d} 
Let $G$ be a locally compact group and $\Ga$ be a generalised $FC^-$-group which acts on $G$ by 
automorphsm. Let $H$ be a closed normal $\Gamma$-invariant subgroup. 
Then $\Gamma$-action on $G$ is distal if and only if $\Gamma$-actions on both 
$H$ and $G/H$ are distal.
\et

\bo 
Let $G$, $H$ and $\Ga$ be as in the hypothesis. Suppose $\Ga$-actions 
on $G/H$ and $H$ are distal. Then it is easy to see that $\Ga$-action on $G$ is 
distal. 

Now we prove the converse. Suppose $\Ga$-action on $G$ is distal and hence $\Ga$-action on
$H$ is also distal. As in the proof of Theorem \ref{tds}, we may assume that $\Ga\subset\Aut(G)$. 
We prove that $\Ga$-action on $G/H$ is distal. By Theorem 3.3 of \cite{RS},
$\Ga$ action on $G/G^0$ is distal and hence equicontinuous (by Proposition \ref{tds}). There exists an open
$\Gamma$ invariant subgroup $L$ in $G$ such that $L=KG^0$, where $K$ is the maximal compact normal
$\Ga$-invariant subgroup of $L$ (cf.\ Lemma \ref{lem}). We know that $G/L$ is discrete, and hence, so is $G/HL$, 
where $HL$ is an open $\Ga$-invariant subgroup. Therefore, it is enough to prove that $\Ga$ acts distally 
on $HL/H$. Since $HL/H$ is isomorphic to $L/(L\cap H)$, without loss of any generality, we may assume
that $G=L=KG^0$ and $K$ is the maximal compact normal $\Ga$-invariant subgroup in $G$. In particular,
$G/K$ is a connected Lie group without any nontrivial  compact normal subgroup.

Here, $HK$ and $K\cap H$ are closed, normal and $\Ga$-invariant subgroups. 
By Theorem 3.1 of \cite{RS} we know that $\Ga$-action is distal on $G/K$, $HK/K$ and on $K/(K\cap H)$; the latter 
is isomorphic to $HK/H$. Hence it is enough to prove that $\Ga$-action is distal on $G/HK$ which is 
isomorphic to $(G/K)/(HK/K)$.  

Replacing $G$ by $G/K$ and $H$ by $HK/K$, we may assume that $G$ is a connected Lie group and 
$H$ is a closed normal (Lie) subgroup and $G$, and hence $H$, has no nontrivial compact normal subgroup.
Let $\cal G$ be the Lie algebra of $G$. Since $\Ga$-action on $G$ is distal, so is the corresponding action of  
$\{{\rm d}\gamma\mid\gamma\in\Ga\}$ on $\cal G$ (cf.~\cite{Ab2}, Theorem 1.1). Equivalently, 
the eigenvalues of ${\rm d}\gamma$ are of absolute value 1, for all $\gamma\in\Ga$ (cf.~\cite{Ab1}, Theorem $1'$). 
Since $H$ is normal and $\Gamma$-invariant, the Lie algebra $\cal H$ of $H^0$ is a Lie 
subalgebra which is an ideal invariant under ${\rm d}\gamma$, for all $\gamma\in\Ga$, and the Lie algebra of 
$G/H$ is isomorphic to $\cal G/\cal H$. Then the eigenvalues of 
${\rm d}\gamma$ on $\cal G/\cal H$ are also of absolute value 1 for all 
$\gamma\in\Gamma$. Hence $\Ga$-acts distally on $G/H$ (cf.~\cite{Ab1}, \cite{Ab2}). 
\hfill{$\Box$}
\eo

\end{section}
\begin{section}{Distality and [MOC]}

In this section we show that distality and [MOC] of  $\Gamma$-action on any locally compact group are equivalent if $\Gamma$  
 is a locally compact, compactly generated abelian (resp.\ Moore) group acting on the group by automorphisms. We first prove a 
 proposition which will be useful in proving Theorem \ref{main}.

\bp\label{close-cpt} Let $G$ and $\Gamma$ be as in Theorem \ref{main}. Suppose
that the $\Gamma$-action on $G$ is distal. 
Given a net $\{\gamma_d\}$ in $\Gamma$, let 
$$M=\{g\in G\mid \{\gamma_d(g)\}_d\mbox{ is relatively compact}\}.$$ 
Then $M$ is a closed $\Ga$-invariant subgroup.
\ep

\bo It is obvious that $M$ is a subgroup and it is $\Ga$-invariant since $\Ga$ is abelian.
Therefore $\ol{M}$ is also a $\Ga$-invariant subgroup. 
If $M$ is trivial, then $M=\ol{M}$. Suppose $M$ is a nontrivial subgroup of $G$. 
Without loss of any generality, we may assume that $G=\ol{M}$, i.e.\ $M$ is dense in $G$.

\medskip
\noindent{\bf Step 1} By Theorem 3.3 of \cite{RS}, $\Ga$-action on $G/G^0$ is distal.
Since $\Ga$ is a compactly generated locally compact abelian group, it 
is a generalised $FC^-$-group. By Proposition \ref{tds}, $\Ga$-action on $G/G^0$ has [MOC] and
$\Ga$-action on $G/G^0$ is equicontinuous.  By Lemma \ref{lem}, there exists an open (resp.\ compact)
$\Gamma$-invariant subgroups $H$ (resp.\  $K$) such that $H=KG^0$, where $K$ is the maximal compact 
normal subgroup of $H$.  Since $H$ is open and $\Ga$-invariant, it is enough to show
that $H\subset M$ and hence, we may assume that $G=H$. Here, since $K$ is a maximal compact normal
$\Ga$-invariant subgroup, $K\subset M$ and $G/K$ is a connected Lie group without
any nontrivial compact subgroup. Moreover, $\Ga$ action on $G/K$ is distal (cf.\ \cite{RS}, Theorem 3.1).
Let $\pi:G\to G/K$ be the natural projection. Since $K$ is compact, 
$\pi(M)=\{gK\in G/K\mid\{\gamma_d(gK)\}_d\mbox{ is relatively compact  in }G/K\}$ 
and $M$ is closed if and only if $\pi(M)$ is closed. Moreover, $\ol{\pi(M)}$ is dense in $G/K$. 
Now, we may replace $G$ by $G/K$ and assume that
$G$ is a connected Lie group without any nontrivial compact normal subgroup and $\Ga\subset\Aut(G)$. 

\medskip
\noindent{\bf Step 2}
Since $G$ has no nontrivial compact central subgroup,  $\Aut(G)$ is almost 
algebraic (as a subgroup of $GL(\cal G)$) (cf.\ \cite{D1}), where $\cal G$ is 
the Lie algebra of $G$.  Let $\Ga'$ be the smallest almost algebraic subgroup containing $\Ga$ in 
$\Aut(G)$. Here $\Ga'$ is a an open subgroup of finite index in the Zariski closure $\tilde{\Ga}$ of 
$\Ga$ in $GL(\cal G)$, hence $\Ga'$ and $\tilde{\Ga}$ have the same connected component of the identity. 
It follows from Corollary 2.5 of \cite{Ab1}, that the unipotent radical $U$ of $\tilde{\Ga}$ is a closed co-compact 
normal subgroup of $\Ga'$ and as in the proof of Theorem 1.1 in \cite{Ab2}. we have that $U$, and hence $\Ga'$, has 
closed orbits in $G$. 

\medskip
\noindent{\bf Step 3} We now prove that  $\{\gamma_d\}_d$ is relatively compact in $\Aut(G)$. 
Suppose $\{\gamma_d\}_d$ is not relatively compact in $\Aut(G)$. Since $\Aut(G)$ is a Lie group,  
there exists a divergent sequence
$\{\gamma'_n\}$ in the set $\{\gamma_d\}$. We know that $\{\gamma'_n(g)\}$ is relatively compact for
all  $g$ in a dense subgroup $M$.  There exists a countable subgroup 
$M_1\subset M$ which is dense in $G$. Passing to a
subsequence if necessary, we may assume that $\{\gamma'_n(g)\}$ converges for all $g\in M_1$. 

Since $G$ is a Lie group without any compact central subgroup of positive dimension, from Step 2,
for every $g\in M_1$, there exists $\gamma_g\in\Ga'$ such that $\{\gamma'_n(g)\}$ converges 
to $\gamma_g(g)$. Then $\gamma_g^{-1}\gamma'_n(g)\to g$ for every $g\in M_1$.  
Let $V$ (resp.\ $W$) be open relatively compact neighbourhoods of the identity $e$ in $G$ (resp. 
zero in $\cal G$) such that the exponential map from $W$ to $V$ is a diffeomorphism
with $\log$ as its inverse. Let $U$ be an open neighbourhood of $e$ in $G$ such that $\ol{U}\subset V$.
Then $(\du\gamma_g^{-1}\du\gamma'_n)(\log g)\to \log g$, 
and hence $\du\gamma'_n(\log g)\to\du\gamma_g(\log g)$ for all $g\in U\cap M_1$. 

In particular, $\{\du\gamma'_n(w)\}$ converges for all $w$ in a dense subset of $\log U\subset W$ in $\cal G$. 
Since $\cal G$ is a vector space and $\log U$ is open, we get that any dense subset of $\log U$ 
generates $\cal G$ and hence $\{\du\gamma'_n\}$ converges in $GL(\cal G)$. Let $\gamma'$ be the 
limit point of $\{\du\gamma'_n\}$ in $GL(\cal G)$; it is a Lie algebra automorphism. Hence 
$\gamma'=\du\gamma$ for some $\gamma\in \Aut(G)$. Then $\gamma'_n\to\gamma$. This is a 
contradiction to the above assumption that $\{\gamma'_n\}$ is divergent. Hence we have
that $\{\gamma_d\}$ is relatively compact in $\Aut(G)$.  This implies that 
$\{\gamma_d(x)\}_d$ is relatively compact for all $x\in G$ and $G=M$, i.e.\ $M$ is closed. \hfill{$\Box$}
\eo

\br\label{rem2} From the above proof it is clear that if $G$ is a connected Lie group without any 
nontrivial compact central subgroup and if $\Ga$ is a subgroup of $\Aut(G)$
acting distally on $G$ and if $\{\gamma_d\}\subset\Ga$ is such that $\{\gamma_d(g)\}_d$ is relatively 
compact for all $g$ in a dense subgroup of $G$, then $\{\gamma_d\}$ is relatively 
compact in $\Aut(G)$.
\er 

\bo $\!\!\!\!\!$ {\bf of Theorem \ref{main}} 
Let $G$ be a locally compact group and let $\Gamma$ be a compactly generated locally compact 
abelian group. Suppose $\Ga$-action on $G$ has [MOC], then we know that $\Ga$-action on $G$ is distal. 

Now suppose that the $\Ga$-action on $G$ is distal. We show that it has [MOC].  
Let $x\in G$ and let $y\in\ol{\Ga(x)}$. We need to show that
$x\in\ol{\Ga(y)}$. We have that $\gamma_d(x)\to y$ for some $\{\gamma_d\}\subset\Ga$. Let 
$$
M=\{g\in G\mid \{\gamma_d(g)\}_d\mbox{ is relatively compact}\}.$$ 
By Proposition \ref{close-cpt}, 
$M$ is a closed $\Ga$-invariant subgroup and $x$, and hence, $y$ belongs to $M$. Without 
loss of any generality we may assume that $M=G$. In view of Theorem \ref{cpt} and Remark \ref{rem1}, 
we can go modulo the maximal compact normal subgroup of $G^0$ which is characteristic in $G$ 
and assume that $G^0$ is a Lie group without any nontrivial compact normal subgroup.  Note that $\Ga$ is
a generalised $FC^-$-group and $\Ga$-action on $G/G^0$ is distal (by Theorem 3.3 of \cite{RS}). Hence from
Proposition \ref{tds}, we get that the action of $\Ga$ on $G/G^0$ is equicontinuous. Let $H_d=K_d\times G^0$ be open 
$\Ga$-invariant subgroups, where $K_d$ are totally disconnected compact $\Ga$-invariant subgroups such that
 $\cap_d K_d=\{e\}$ in $G$ (cf.\ Lemma \ref{lem}). Then passing to a subnet if necessary, we may assume that 
 $\gamma_d(x)=yk_dg_d=yg_dk_d$, 
where $k_d\in K_d$ and $g_d\in G^0$, $k_d\to e$, $g_d\to e$. In particular, we get that 
$\gamma_d^{-1}(y)=x\gamma_d^{-1}(k_d^{-1})\gamma_d^{-1}(g_d^{-1})$.
We know that $\{\gamma_d|_{G^0}\}$ is relatively compact, (see Remark \ref{rem2}). Let $\ga$ be a limit
point of $\{\gamma_d|_{G^0}\}$ in $\Aut(G^0)$. Then $\gamma^{-1}$ is a limit point of $\{\gamma_d^{-1}|_{G^0}\}$ in
$\Aut(G^0)$.  Therefore, passing to a subnet if necessary, we get that 
$$
\gamma_d^{-1}(g_d^{-1})\to\gamma^{-1}(e)=e\ \ \mbox{and}\ \ \gamma_d^{-1}(y)=xk'_d\gamma_d^{-1}(g_d^{-1})\to x$$
where $k'_d=\gamma_d^{-1}(k_d^{-1})\in K_d$ and $k'_d\to e$ as $K_d$ are $\Ga$-invariant and $\cap_d K_d=\{e\}$. 
In particular, $x\in\ol{\Ga(y)}$. 

Since the above is true for any $x\in G$ and any $y\in\ol{\Ga(x)}$, closure of any orbit is a
minimal closed $\Ga$-invariant set, i.e.\ $\Ga$-action on $G$ has [MOC]. \hfill{$\Box$}
\eo

A locally compact group $G$ is said to be a {\it central} group or a $Z$-group if $G/Z(G)$ is compact, 
where $Z(G)$ is the center of $G$. It is said to be a {\it Moore} group if all its irreducible unitary representations
are finite dimensional. All abelain groups and all compact groups are $Z$ groups and $Z$-groups are also 
Moore groups.  A Moore group has normal subgroup $H$ of finte index such that $\ol{[H,H]}$ is compact 
(cf.\ \cite{Rob}). It is easy to see from this, that any Moore group $G$ is $FC^-$ nilpotent as 
$G_0=G$, $G_1=H$, $G_2=\ol{[H,H]}$ and $G_3=\{e\}$. Since $G_0/G_1$ is finite, and 
$G_1/G_2$ is abelian and $G_2/G_3$ is compact, we have that the conjugacy action of $G$  on 
$G_i/G_{i+1}$ has relatively compact orbits for all $i=0,1,2$.  Hence any compactly generated 
Moore group has polynomial growth and it is a generalised $FC^-$-group (cf.\ \cite{L2}, 
Theorem 1, Lemma 1). 

\bc\label{mainc} Let $G$ be a locally compact group and let $\Ga$ is a compactly generated Moore 
group acting on $G$ by auotomorphisms. 
Then $\Ga$-action on $G$ is distal if and only if it has {\rm [MOC]}. 
\ec

The proof of the above corollary is essentially the same as that of Theorem \ref{main}. As $\Ga$ is a Moore 
group, it has a closed normal subgroup $\Ga_1$ of finite index whose commutator group is relatively compact.
(cf.\ \cite{Rob}, Theorem 1). Then by Lemma 4.1 of \cite{Mo}, it is enough to show that $\Ga_1$-action on
$G$ has [MOC]. Without loss of any generality, we may assume that 
$[\Ga,\Ga]$ is relatively compact and hence it is easy to see that the 
group $M$ defined in the above proof is $\Ga$-invariant. We will not repeat the proof here. 

\br\label{rem3} 1. From above, it is obvious that Theorem \ref{main} holds for any
compactly generated locally compact group $\Ga$ such that its commutator subgroup is relatively compact. 
Moreover from Lemma 4.1 in \cite{Mo} we know that the action of a group $\Gamma$ on $G$ has 
[MOC] if  the action of any co-compact subgroup of $\Ga$ on $G$ has [MOC]. Hence Theorems
\ref{main} and \ref{normal} hold for compact extensions of such a group $\Ga$ mentioned above, 
and in particular, for compact extensions of compactly generated abelian, or more generally, of Moore 
groups.  
\er

We conjecture that Theorem \ref{main} holds for any generalised $FC^-$-groups. 
It already holds for the actions of such a group on totally disconnected groups, compact groups 
and connected groups.

\end{section}

\begin{acknowledgement}
The author would like to thank H.\ Abels, S.G.\ Dani, Y.\ Guivarc'h and C.R.E.\ Raja for fruitful
discussions. The author would also like to thank C.R.E.\ Raja for comments on a preliminary version
which led to improvement in the presentation of the manuscript. 
\end{acknowledgement}

\bigskip
\advance\baselineskip by 1pt
\noindent Riddhi Shah\\
School of Physical Sciences(SPS)\\
Jawaharlal Nehru University(JNU)\\
New Delhi 110 067, India\\
rshah@mail.jnu.ac.in\\
riddhi.kausti@gmail.com

\end{document}